\documentclass[11pt,a4paper]{article}
\usepackage{epsfig}
\usepackage{times}
\usepackage{amssymb} 
\newenvironment{keywords}{ \noindent {\small\bf Key Words}:}{ }

\newcommand{\be}{\begin{equation}}
\newcommand{\ee}{\end{equation}}

\newcommand{\bdm}{\begin{displaymath}}
\newcommand{\edm}{\end{displaymath}}

\newcommand{\bd}{\begin{description}}
\newcommand{\ed}{\end{description}}

\newtheorem{theorem}{Theorem}
\newtheorem{corollary}{Corollary}

\begin{document}

\title{Local tuning and partition strategies for \\ diagonal GO
methods \thanks{This research was partially supported by the
Russian Fund of Basic Research through grant number
01--01--00587}}

\author{Dmitri E. Kvasov\\
University of Nizhni Novgorod,   Nizhni
Novgorod, Russia  \\
 Clara Pizzuti\\
Istituto per la Sistemistica e l' Informatica, C.N.R., c/o
D.E.I.S. \\ Universit\`{a} della Ca\-lab\-ria, 87036 Rende (CS),
Italy \\
 Yaroslav D. Sergeyev\thanks{Corresponding author
 \textit{yaro@si.deis.unical.it}}\\
University of Nizhni Novgorod,   Nizhni
Novgorod, Russia  and\\
 Istituto per la Sistemistica e l'
Informatica, C.N.R., c/o D.E.I.S. \\ Universit\`{a} della
Ca\-lab\-ria, 87036 Rende (CS), Italy, }
 \date{}
%

\maketitle

\begin{abstract}
In this paper, global  optimization (GO) Lipschitz problems are
considered where the multi-dimensional multiextremal objective
function is determined over a hyperinterval. An efficient
one-dimensional GO method using local tuning on the behavior of
the objective function is generalized to the multi-dimensional
case by the diagonal approach using two partition strategies.
Global convergence conditions are established for the obtained
diagonal geometric methods. Results of a wide numerical comparison
show a strong acceleration reached by the new methods working with
estimates of the local Lipschitz constants over different
subregions of the search domain in comparison with the traditional
approach.
\end{abstract}



\section{Introduction} \label{intro}

In \cite{Pinter (1983),Pinter (1986),Pinter (1996)} diagonal
global optimization algorithms have been introduced for solving
multi-dimensional Lipschitz global optimization (GO) problems with
box constraints. In its general form such a problem can be stated
as minimization of a multiextremal function satisfying the
Lipschitz condition with a constant $0<L<\infty$ over a
hyperinterval, i.e., finding the value $f^*$ and points $x^*$
such that
 \be
 f^*=f(x^*)=\min_{x\in D } f(x), \label{1}
 \ee
 where
 \be
 |f(x')-f(x'')| \le L \parallel x'- x''\parallel, \hspace{5mm}
x', x'' \in D \subset \mathbb{R}^n,        \label{1.1}
 \ee
 \be D=[a,b\,]= \{x \in \mathbb{R}^n\,:\, a \leq x \leq b\},
 \hspace{5mm} a \leq b,\ \ a,b \in \mathbb{R}^n.   \label{2}
 \ee

Such problems very often can be faced in real-life applications
(for example, in data classification, nonlinear approximation,
globally optimized calibration of complex system models etc.). A
number of such problems solved by the diagonal methods can be
found in \cite{Pinter (1996)}.

The diagonal approach is a simple and powerful tool for extending
one-dimensional global optimization methods to the
multi-dimensional case. The main idea is to describe the behavior
of the objective function $f(x)$ over a hyperinterval (we shall
also use the term {\it cell} or simply {\it interval})
$D_i=[a_i,b_i]$ by information obtained from evaluating $f(x)$ at
the vertices $a_i$, $b_i$ being the ends of the main diagonal
defining the interval $D_i$. During every $(l+1)$-th iteration to
each subinterval $D_i \subset D$ generated in the course of the
previous $l$ iterations a {\it characteristic}
$R_i=R(a_i,b_i,f(a_i),f(b_i))$ is associated in such a way that
$R_i$ tends to be higher if $D_i$ contains lower values of $f(x)$.
Then, among all subintervals created so far within $D$, an
interval $D_t$ with the maximal characteristic is chosen for
further subdivision. It is subdivided in $p$ subcells and $f(x)$
is evaluated at the vertices $a_j$, $b_j$ of all the intervals
$D_j, 1 \le j \le p$. The process is repeated until satisfaction
of a stopping rule.

The diagonal method proposed in \cite{Pinter (1983),Pinter
(1986),Pinter (1996)} and extending the univariate algorithm from
\cite{Pijavskii (1972)} uses a global estimate of the Lipschitz
constant $L$ in its work. GO algorithms using in their work the
global Lipschitz constant $L$ (or its estimates) do not take into
account local information about behavior of the objective function
over every small subregion of $D$. In fact, it is supposed in such
algorithms (see  \cite{Horst and Pardalos (1995)}) that $f(x)$ has
the  same constant $L$ over every subdomain of $D$ without paying
any attention to situations where $f(x)$ has a very low  local
Lipschitz constant over the subdomain under consideration. It has
been shown  for a number of global optimization algorithms (see
\cite{Sergeyev (1995a),Sergeyev (1995b),Strongin&Sergeyev (2000)})
that using local information for estimating  local Lipschitz
constants can accelerate the global search significantly.
Importance of such information in the diagonal approach context
has been highlighted in \cite{Pinter (1996)}. Of course, the local
data must be in an appropriate way balanced with the global
information about the objective function otherwise the global
solution can be lost \cite{Stephens (1998)}.

In this paper a new diagonal algorithm generalizing an efficient
deterministic one-dimensional GO method using local tuning on the
behavior of the objective function (see \cite{Sergeyev (1995b)})
is extended to the multi-dimensional case by the diagonal approach
using two partition strategies widely used in literature
\cite{Gergel (1997),Meewella and Mayne  (1989),Pinter
(1983),Pinter (1986),Pinter (1996)}:
 \bd
 \item
 -- Bisection,  where $p=2$ and the interval $D_t$ is subdivided in two
 subintervals by a hyperplane orthogonal to the longest edge of $D_t$;
 \item
 -- Partition $2^n$, where $p=2^n$ and
$D_t$ is partitioned into $2^n$ new subintervals  generated by the
intersection of the boundary of $D_t$ and the hyperplanes that
contain a point $x^{l+1}$ belonging to the main diagonal of $D_t$
and are parallel to the boundary hypersurfaces of $D_t$.
 \ed

The new method uses a local information about the objective
function over the {\it whole} search region $D$ {\it during}  the
global search in contrast with techniques which do it only in a
{\it neighborhood} of local minima {\it after} stopping their
global procedures (see e.g. \cite{Horst and Pardalos (1995)}).
Global convergence conditions are established for the new method.
Results of a wide numerical comparison show a strong acceleration
reached by the new method working with estimates of the local
Lipschitz constants over different subregions of the search
domain in comparison with the traditional approach using global
estimates of $L$.

\section{The new algorithm  with local tuning} \label{sec:2}

In this section the New Diagonal Algorithm  with Local tuning
(NDAL) is described.

The method starts by setting  the number of iterations, $l$, and
the number of generated intervals, $m=m(l)$, equal to 1. The first
two {\it trials} (evaluations of the objective function) are
executed at the points $x_0 = a$,  $x_1 = b$ from (\ref{2}). The
results of trials are indicated as $z_0 = f(x_0)$, $z_1 =
f(x_1)$, and the initial number $k=k(l)$ of trial points generated
by the algorithm is taken equal to 2. The initial estimate of the
global optimum is taken as $ z^*_1 = \min\{z_0,z_1\}$. The
estimate $\lambda_1$ of the local Lipschitz constant over the
initial interval $D_1=D=[a,b\,]$ (in this case, of course, the
local estimate coincides with the global one) is calculated as
follows
 \bdm
  \lambda_1 = \frac {\mid f(a)-f(b) \mid}{\parallel a-b
  \parallel}.
 \edm

Suppose now that $l\ge1$ iterations of the method have already
been executed. The iteration $l+1$ consists of the following
steps.
 \bd
\item {\bf Step 1}.    For each
interval $D_i = [a_i,b_i]$, $1\leq i\leq m(l)$, calculate its
characteristic
 \be
  R_i=0.5(K_i\parallel a_i-b_i\parallel
  -f(a_i)-f(b_i))    \label{13}
 \ee
where
 \be
  K_i= K_i(l)=(r+\frac Cl) \max \{ \lambda _i, \gamma _i, \xi \}, \label{8}
 \ee
the values  $r>1$, $\xi> 0$, and $C>0$ are parameters of the
method, $\lambda _i$ is the estimate of the local Lipschitz
constant over the interval $D_i$ calculated at the moment of
creation of $D_i$, and
 \be
  \gamma _{i}= \mu \frac{\parallel a_i-b_i\parallel}{d^{\max}}. \label{10}
 \ee
The values  $\mu$ and $d^{max}$ are evaluated as follows
 \be
   \mu = \max_{1\leq i \leq m(l)}  \lambda _i, \label{11}
 \ee
 \be
  d^{\max}= \max_{1\leq i \leq m(l)} \parallel a_i-b_i \parallel. \label{12}
 \ee
\item
{\bf Step 2}. Among all the intervals $D_i$ choose an interval
$D_t$ such that
 \be
 R_t = \max_{1\leq i \leq m(l)}R_i. \label{4}
\ee
\item
{\bf Step 3}.   If
 \bdm
  \parallel a_t-b_t\parallel > \varepsilon \parallel a-b \parallel,
 \edm
where $a$ and $b$ are from (\ref{2}) and $t$ is from (\ref{4}),
then go to Step~4, otherwise take the value
 \bdm
  z^*_l = \min_{1 \leq i \leq k(l)}f(x_i)
 \edm
(where $x_i$, $1 \leq i \leq k(l)$, are the trial points generated
by the algorithm in the course of the previous $l$ iterations) as
an estimate of the global optimum of the
problem~(\ref{1})~--~(\ref{2}) and~{\bf Stop}.
\item
{\bf Step 4}.    Choose the new point $x^{l+1}$ belonging to the
main diagonal (the diagonal joining the vertices $a_t$ and $b_t$)
of the subinterval $D_t$, where $t$ is from (\ref{4}), as follows
(see \cite{Pinter (1983),Pinter (1986),Pinter (1996)}):
 \be
x^{l+1}=\frac{a_t+b_t}{2} - \frac{f(b_t)-f(a_t)}{2\hat{K}} \times
\frac{b_t-a_t}{\parallel a_t-b_t\parallel}. \label{14}
 \ee
Here
 \be
  \hat{K} = \hat{K}(l) = (4+\frac{C}{l})\max \{ \mu, \xi \} ,  \label{6}
 \ee
where $\xi$ is from (\ref{8}) and $\mu$ is from (\ref{11}).
\item
{\bf Step 5}.   Subdivide the interval $D_t$ into $p$ new
subintervals by Bisection strategy or by Partition $2^n$.
 \item
{\bf Step 6}.   Denote by $x_i$, $i=1,\ldots ,s$, the vertices
of the new $p$ subintervals generated during Step 5 where $f(x)$
must be evaluated.

--In the case of Bisection strategy it is necessary to evaluate
$f(x)$ at two vertices, $s=2$ (the points $a_t$ and $b_t$ come
from the subdivided interval $D_t$  and $f(x)$ has already been
evaluated at its vertices during the previous iterations).

--In the case of Partition $2^n$, the number $s=2 \times 2^n -3$
because the new $2^n$ subintervals are identified by their two
vertices, $x^{l+1}$ is common to two intervals, and $f(a_t)$ and
$f(b_t)$ of the subdivided interval $D_t$ have already been
evaluated.
\item
{\bf Step 7}.   For all the new intervals $D_i$, $1 \leq i \leq
p$, get an estimate of the local Lipschitz constant as
 \be
  \lambda _i= \max \{\frac {\mid f(a_t)-f(b_t) \mid}
     {\parallel a_t-b_t\parallel}, \max_{1 \leq j \leq p} \frac
     {\mid f(a_j)-f(b_j) \mid} {\parallel a_j-b_j\parallel}\}.
     \label{9}
 \ee

Set $l:=l+1$, $m:=m+p-1$, $k:=k+s$, and go to Step~1.

  \ed

Let us give a few comments on the introduced method. The key idea
of the algorithm is estimating local Lipschitz constants by
balancing local and global data. In contrast with the traditional
approach (see \cite{Pinter (1983),Pinter (1986)}) where the global
estimate $\hat{K}$ of the Lipschitz constant $L$ from (\ref{1.1})
is used in the form (\ref{6}), the local estimate $K_i$ from
(\ref{8}) is the result of the balance between the local and the
global information represented by the values $\lambda _{i}$ and
$\gamma _{i}$, respectively. When the subinterval $D_{i}$ has a
small main diagonal (in comparison with the current maximal
diagonal $d^{\max}$ over all subintervals in $D$)  then (see
(\ref{10})--(\ref{12})), $\gamma _{i}$ is small too and   the
local information represented by $\lambda _{i}$ has a decisive
influence (see (\ref{8})) on $K_i$. When the interval $D_{i}$ is
very wide (its diagonal $\parallel a_i-b_i\parallel$ is close to
$d^{\max}$), the local information is not reliable and the global
information (see (\ref{10})) represented by $\gamma _{i}$ is used.

The values $r$, $C$, and $\xi $ influence $K_i$ as global
parameters. By increasing $r$ and $C$ we augment  reliability of
the method over the whole region $D$. The parameter $\xi  > 0$ is
a small number allowing the NDAL to work also when $f(x_i)=const$
for all trial points $x_i$. The importance of the parameter $\xi$
for the correct  work of the method can be seen from
(\ref{13})~--~(\ref{8}) and (\ref{14})~--~(\ref{6}). If $\gamma
_{i}< \xi$ and $\lambda _{i} < \xi$ it follows
 \bdm
  K_i(l) = \hat{K}(l) = (r+\frac Cl)  \xi.
 \edm
Of course, this case is degenerate for the method.

The introduced algorithm belongs to the  class of diagonally
extended geometric algorithms and also to more general classes of
{\it adaptive partition} and {\it divide the best} algorithms
(see \cite{Pinter (1992),Pinter (1996)} and \cite{Sergeyev
(1999)}, respectively). Let us study the convergence properties of
the infinite ($\varepsilon=0$ in the stopping rule) sequence
$\{y^{k}\}$ of trial points generated by the NDAL during
minimization of the function $f(x)$ from (\ref{1})--(\ref{2}).
Hereinafter we shall designate by $Y'$ the set of limit points of
the sequence $\{y^{k}\}$.

\begin{theorem}
Let $y'$ be a limit point  of the sequence $\{y^{k}\}$ then, for
all trial points $y^{k}\in \{y^{k}\}$, it follows   $f(y^{k})\geq
f(y')$. If there exists another limit point $y'' \in Y'$ then
$f(y') = f(y'')$. \label{t2}
\end{theorem}
{\bf Proof.} This result can be obtained as a particular case of
the general convergence study from \cite{Sergeyev (1999)} and its
proof  is so omitted. $\Box$

The next theorem presents sufficient global convergence
conditions for the NDAL.
\begin{theorem}
Let  there exist an iteration  number $l^*$ such that for a cell
$D_{j}$, $j=j(l)$,  containing a global minimizer $x^{*}$ of
$f(x)$ during the $l$-th  iteration of the NDAL the following
inequality takes place
 \be K_{j}(l) \ge  2H_{j},\hspace {1cm} l > l^*, \label{21}
\ee
 where
  \be H_{j}=
\max \{ \frac{f(a_{j})- f(x^{*})} {\parallel x^{*} -
a_{j}\parallel}, \frac {f(b_{j})- f(x^{*})} {\parallel b_{j}-
x^{*}\parallel} \}. \label{22}
 \ee
Then, $x^{*}$ is a limit point of the trial sequence $\{y^{k}\}$
generated by the NDAL. \label{t3}
\end{theorem}
{\bf Proof.} We start the proof by showing  that the estimates
$K_{i}(l)$ of the local Lipschitz constants $L_{i}$ from
(\ref{8}) are bounded values. In fact, since the global Lipschitz
constant $L < \infty$ and the constants $r > 1$, $C > 0$, and $\xi
> 0$, it follows
 \be
  0 <
r\xi < K_{i}(l) \leq (r+C) \max \{ L, \xi\} < \infty,
\hspace{1cm} l \geq 1. \label{23}
 \ee
Suppose, that there exists a limit point $y' \neq x^{*}$ of the
trial sequence $ \{ y^{k} \}$. Taking into consideration
(\ref{13}), (\ref{14}), (\ref{6}), and (\ref{23}) we can conclude
for an interval $D_i$, $i =i(l)$, containing $y'$ during the
$l$-th iteration of the NDAL, that
 \be
    \lim_{l \rightarrow \infty} R_i(l)=-f(y').  \label{24}
 \ee
Consider now the cell $D_{j}$, $j=j(l)$, such that the global
minimizer $x^{*} \in D_{j}$ and suppose that $x^{*}$ is not a
limit point of $\{y^{k} \}$. This signifies that there exists an
iteration number $q$ such that for all $l \geq q$
 \bdm
  x^{l+1} \notin D_{j}, \hspace{5mm} j=j(l).
 \edm
Estimate now the characteristic $R_j(l)$, $l \geq q$, of the
interval $D_{j}$. It follows from (\ref{22}) and the fact of
$x^{*} \in D_{j}$ that
 \bdm
  f(a_{j}) - f(x^{*}) \leq H_{j}\parallel a_{j}-x^{*}\parallel \leq
  H_{j}\parallel a_{j}- b_{j}\parallel ,
 \edm
 \bdm
  f(b_{j})-f(x^{*}) \leq H_{j}\parallel b_{j}-x^{*}\parallel \leq
  H_{j}\parallel a_{j}- b_{j}\parallel .
 \edm
Then, by summarizing these inequalities  we obtain
 \bdm
  f(a_{j})+f(b_{j}) \leq 2f(x^{*}) + 2H_{j}\parallel a_{j}-
  b_{j}\parallel.
 \edm
From this inequality and (\ref{21}), (\ref{22}) we can deduce for
all iteration numbers $l > l^*$ that
 \bdm
  R_j(l)=0.5(K_{j}\parallel a_{j}- b_{j}\parallel -
  f(a_{j})-f(b_{j})) \geq
 \edm
 \bdm
  0.5(K_{j}\parallel a_{j}- b_{j}\parallel
  -2f(x^{*})-2H_{j}\parallel a_{j}- b_{j}\parallel) =
 \edm
 \be
  0.5\parallel a_{j}- b_{j}\parallel(K_{j} -2H_{j})-f(x^{*}) \geq -f(x^{*}).
  \label{25}
 \ee
Since $x^{*}$ is a global minimizer, it follows from (\ref{24})
and (\ref{25}) that an iteration number $q^{*} > \max \{l^{*},q
\}$ will exist such that
 \bdm
  R_j(q^{*}) \geq R_i(q^{*}).
 \edm
But this means that during the $q^{*}$-th iteration, trials will
be executed at the cell $D_j$.  Thus, our assumption that $x^{*}$
is not a limit point of $\{y^{k} \}$ is not true and theorem has
been proved.  $\Box$

Let us denote the set of global minimizers of the problem
(\ref{1})--(\ref{2}) as $X^{*}$. Then the following corollary
ensures the inclusion $Y'\subseteq X^{*}$.

\begin{corollary}
Given the conditions of  Theorem~\ref{t3}, all limit points of the
sequence $\{y^{k} \}$ are global minimizers of $f(x)$,
$Y'\subseteq X^{*}$. \label{c3}
\end{corollary}
{\bf Proof.} The corollary follows immediately from
Theorems~\ref{t2} and~\ref{t3}. $\Box$

The  sets $Y'$ and $ X^{*}$ coincide if  conditions established
by Corollary~\ref{c4} are fulfilled.
\begin{corollary}
If condition (\ref{21}) is fulfilled for all points  $x^{*} \in
X^{*}$, then the set of limit points of $\{y^{k} \}$ coincides
with the set of global minimizers of the objective function
$f(x)$, i.e. $Y' = X^{*}$. \label{c4}
\end{corollary}
{\bf Proof.} Again, the corollary is a straightforward consequence
of Theorems~\ref{t2} and~\ref{t3}. $\Box$

\section{Numerical comparison} \label{sec:3}

The goal of this section is dual: first, to show advantages of the
local tuning in comparison to the traditional approach using
global estimates of the Lipschitz constant; second, to establish
which of two partitioning strategies, Bisection or Partition
$2^n$, works better.

\newpage

Thus,  four methods are compared:
 \bd
 \item
 -- the traditional method with Partition $2^n$ and the global estimate;
 \item
  -- the traditional method with Bisection and the global estimate;
 \item
  -- the new algorithm using local tuning and Partition $2^n$;
  \item
  -- the new algorithm using local tuning and Bisection.
 \ed

\begin{table}[t]
\caption{Test problems} \fontsize{7}{8,4} \label{lstproblems}
\begin{tabular}{cccc}
\hline\noalign{\smallskip}
$N^o$&Formula&Domain&Source\\
\noalign{\smallskip}\hline
1&$0.25x^4_1 - 0.5x^2_1 + 0.1x_1 + 0.5x^2_2$ &$[-10,10]^2$&\cite{Lucidi (1989)}\\
\hline 2&$(4-2.1x^2_1 +x^4_1/3)x^2_1 + x_1 x_2 + (-4 +
4x^2_2)x^2_2$&$[-2.5,2.5]\times$&\cite{Torn and Zilinskas
(1989)}\\
        & &$[-1.5,1.5]$& \\
\hline 3&$2x^2_1 - 1.05x^4_1 + x^6_1/6 + x_1x_2
+x^2_2$&$[-5,5]^2$&\cite{Dixon and Szego (1975)}\\
\hline 4&$(x_2 - 5.1x^2_1/(4\pi ^2) + 5x_1/\pi - 6)^2 + 10(1-
1/(8\pi))\cos x_1
+10$&$[-5,10] \times$&\cite{Branin (1972)}\\
        & &$[\,0,15]$& \\
\hline 5&$( 1 - 2x_2 + 0.05\sin (4{\pi}x_2) - x_1)^2 + ( x_2 -
0.5\sin (2{\pi}x_1)
)^2$&$[-10,10]^2$&\cite{Dixon and Szego (1975)}\\
\hline 6&$[1 +
(x_1+x_2+1)^2(19-14x_1+3x^2_1-14x_2+6x_1x_2+3x^2_2)]\times$
&$[-2,2]^2$&\cite{Goldstein and Price (1971)}\\
       &$[30+(2x_1-3x_2)^2(18-32x_1+12x^2_1+48x_2-36x_1x_2+27x^2_2)]$& &\\
\hline 7&$\sum_{i=1}^{5} i \cos ((i+1)x_1+i) \sum_{j=1}^{5} j \cos
((j+1)x_2+j)$&$[-10,10]^2$&\cite{Lucidi (1989)}\\
\hline 8&$\sum_{i=1}^{5} i \cos ((i+1)x_1+i) \sum_{j=1}^{5} j
\cos ((j+1)x_2+j)
+$&$[-10,10]^2$&\cite{Lucidi (1989)}\\
        &$(x_1+1.42513)^2 + (x_2 + 0.80032)^2$&&\\
\hline 9&$100(x_2-x^2_1)^2+(x_1-1)^2$&$[-2,8]^2$&\cite{Dixon and
Szego
(1975)}\\
\hline
10&$(x^2_1+x_2-11)^2+(x_1+x^2_2-7)^2$&$[-6,6]^2$&\cite{Himmelblau
(1972)}\\ \hline 11&$-4x_1x_2\sin
(4{\pi}x_2)$&$[\ 0,1]^2$&\cite{Mladineo (1986)}\\
\hline 12&$-\sin
(2x_1+1)-2\sin (3x_2+2)$&$[\ 0,1]^2$&\cite{Mladineo (1986)}\\
\hline
13&$(x_1-2)^2+(x_2-1)^2-0.04/(0.25x^2_1+x^2_2-1)+5(x_1-2x_2+1)^2$&$[\
1,2]^2$&\cite{Schittkowski (1987)}\\
\hline 14&$-\mid\!\sin (x_1)\sin (2x_2)\!\mid+0.01(x_1x_2+(x_1-\pi
)^2+3(x_2-\pi )^2)$&$[\ 0,2\pi ]^2$&\cite{Sergeyev et al. (2001)}\\
\hline 15&$(\pi/n)\{10 \sin ^2(\pi y_1) + \sum^{n-1}_{i=1}[(y_i
-1)^2(1 + 10 \sin ^2(\pi y_{i+1}))] + $ &$[-10,10]^n$&\cite{Lucidi (1989)}\\
      &$(y_n-1)^2\},\ where\  y_i = 1+(1/4)(x_i-1)$, \hspace{2mm} $i=1, \ldots , n$& &\\
\hline 16&$0.1\{\sin ^2(3 \pi x_1) + \sum^{n-1}_{i=1}[(x_i
-1)^2(1 +  \sin ^2(3 \pi
x_{i+1}))]\} + $ & $[-10,10]^n$ &\cite{Lucidi (1989)}\\
      &$0.1(x_n -1)^2 [1+ \sin ^2(2 \pi x_n)\;]$& & \\
\hline
17&$-\sum^{4}_{i=1}c_i\exp\,(-\sum^{3}_{j=1}\alpha_{ij}(x_j-p_{ij})^2)$&$[\
0,1]^3$&\cite{Hartman (1973)}\\
\hline 18&$100[x_3-0.25(x_1+x_2)^2]^2+(1-x_1)^2+(1-x_2)^2$&$[\
0,1]^3$&\cite{Schittkowski (1987)}\\
\hline 19&$(x^2_1-2x^2_2+x^2_3)\sin (x_1)\sin (x_2)\sin
(x_3)$&$[-1,1]^3$&\cite{Mladineo (1992)}\\
\hline
20&$\sum^{3}_{i=1}[(x_1-x^2_i)^2+(x_i-1)^2]$&$[-10,10]^3$&\cite{Walster
et al. (1985)}\\
\hline
\end{tabular}
\end{table}

The list of problems used in the experiments is shown in
Table~\ref{lstproblems}, where the following quantities are
specified:
 \begin{description}
 \item $N^o$ : problem number;
 \item $Formula$ : formula of the test function;
 \item $Domain$ : feasible region of the test function;
 \item $Source$ : bibliographic reference.
 \end{description}

Problems 1--14 are two\--di\-men\-sio\-nal, problems 17--20 are
three\--di\-men\-sio\-nal, and problems 15--16 are of arbitrary
dimension $n>1$ ($n=2$ and $n=3$ have been used).

To show the influence of the parameter $r$ on the search
characteristics, the experiments for the two-dimensional case
have been realized for two different values of the parameter $r$
in all the methods: $r=1.1$ and $r=1.3$. The value $C=10$ was
taken in all the two-dimensional experiments. We have executed
these experiments with the accuracy $\varepsilon=0.01$ in the
stopping rule.

\begin{table}[t]
\caption{Results of numerical experiments with two-dimensional
functions for $r=1.1$} \label{tab2r11}
\begin{center}
\begin{tabular}{@{\extracolsep{\fill}}|c|r|r|r|r|}\hline
Problem& \multicolumn{2}{c|}{Global Estimate} & \multicolumn{2}{c|}{Local Tuning} \\
\cline{2-5}Number&Partition $2^n$ &Bisection &Partition $2^n$ &Bisection\\
\hline 1 &12412 &8950  &4742  &3508 \\
\hline 2 &8037  &2670  &2947  &1354 \\
\hline 3 &19427 &20392 &14832 &14244 \\
\hline 4 &4687  &2762 &1332  &998 \\
\hline 5 &4187  &2818  &807   &602 \\
\hline 6 &20522 &17732 &14572  &10924 \\
\hline 7 &6837  &4766  &5532  &3936 \\
\hline 8 &4057  &3922  &2822  &3372 \\
\hline 9 &16187 &16446 &10307  &7328 \\
\hline 10 &6267 &4384  &1797  &1286 \\
\hline 11 &312  &256   &272   &146 \\
\hline 12 &292  &200   &167   &96 \\
\hline 13 &1827 &2002  &282   &238 \\
\hline 14 &1127 &96$^*$   &592   &186 \\
\hline 15 &4857 &2736  &2237  &1336 \\
\hline 16 &1627 &532   &492   &118 \\
\hline Average &7041.36 &5666.50 &3983.25 &3104.50 \\
\hline
\end{tabular}
\end{center}
\end{table}

\begin{table}[t]
\caption{Results of numerical experiments with two-dimensional
functions for $r=1.3$}  \label{tab2r13}
 \begin{center}
\begin{tabular}{@{\extracolsep{\fill}}|c|r|r|r|r|}\hline
Problem&\multicolumn{2}{c|}{Global Estimate} &\multicolumn{2}{c|}{Local Tuning} \\
\cline{2-5}Number&Partition $2^n$ &Bisection &Partition $2^n$ &Bisection\\
\hline 1 &13987 &9874 &7012  &5620 \\
\hline 2 &9862  &4774  &3357  &2072 \\
\hline 3 &20057 &21608 &16802 &16754 \\
\hline 4 &5812  &3728 &2332  &1190 \\
\hline 5 &4817  &3180  &1402  &650 \\
\hline 6 &21922 &22424 &17812  &12622 \\
\hline 7 &7267  &7374  &6422 &5128 \\
\hline 8 &5467  &4504  &3717  &3938\\
\hline 9 &16752 &17378 &10852  &8250 \\
\hline 10 &8852 &6820  &3432  &1858 \\
\hline 11 &417  &324   &362   &174 \\
\hline 12 &347  &232   &177   &114 \\
\hline 13 &2102 &2306  &307   &284 \\
\hline 14 &1297 &800   &747   &360 \\
\hline 15 &7167 &3880  &3137  &1740 \\
\hline 16 &1852 &778   &612   &162 \\
\hline Average &7998.56 &6874.00 &4905.13 &3807.25 \\
\hline
\end{tabular}
 \end{center}
\end{table}

The numbers of function evaluations executed by the methods before
satisfaction of the stopping rule for the two-dimensional case
are reported in Tables~\ref{tab2r11} and~\ref{tab2r13}. Global
optima have been located in all the experiments. For Problem~14
and the method with the global estimate of the Lipschitz constant
and Bisection strategy the value $r=1.1$ was too small: the
method has not located the global minimizer in this case. The
sufficient value of the reliability parameter $r$ for finding the
global minimizer for Problem~14 is~$r=1.3$.

In Table~\ref{tab3r12} the experimental results for
three-dimensional test functions are shown. The following
parameters have been chosen in all the experiments: $r=1.2$,
$C=100$. The search accuracy $\varepsilon = 0.02$ has been used.

\begin{table}
\caption{Results of numerical experiments with three-dimensional
functions for $r=1.2$} \label{tab3r12}
 \begin{center}
\begin{tabular}{@{\extracolsep{\fill}}|c|r|r|r|r|}
\hline Problem&\multicolumn{2}{c|}{Global
Estimate}&\multicolumn{2}{c|}{Local Tuning} \\
\cline{2-5}
Number&Partition $2^n$&Bisection&Partition $2^n$&Bisection\\
\hline 15 &173513 &43780 &98412 &12060  \\
\hline 16 &26938  &3732  &12625  &1032   \\
\hline 17 &6879   &1810  &4825  &1020  \\
\hline 18 &83475  &27760 &15862 &3470  \\
\hline 19 &8556   &2040  &7568  &1358  \\
\hline 20 &122436 &74254 &59646 &21756 \\
\hline Average &70299.50 &25562.67 &33156.33 &6782.67 \\
\hline
\end{tabular}
\end{center}
\end{table}

\begin{figure}
\centerline{\psfig{file=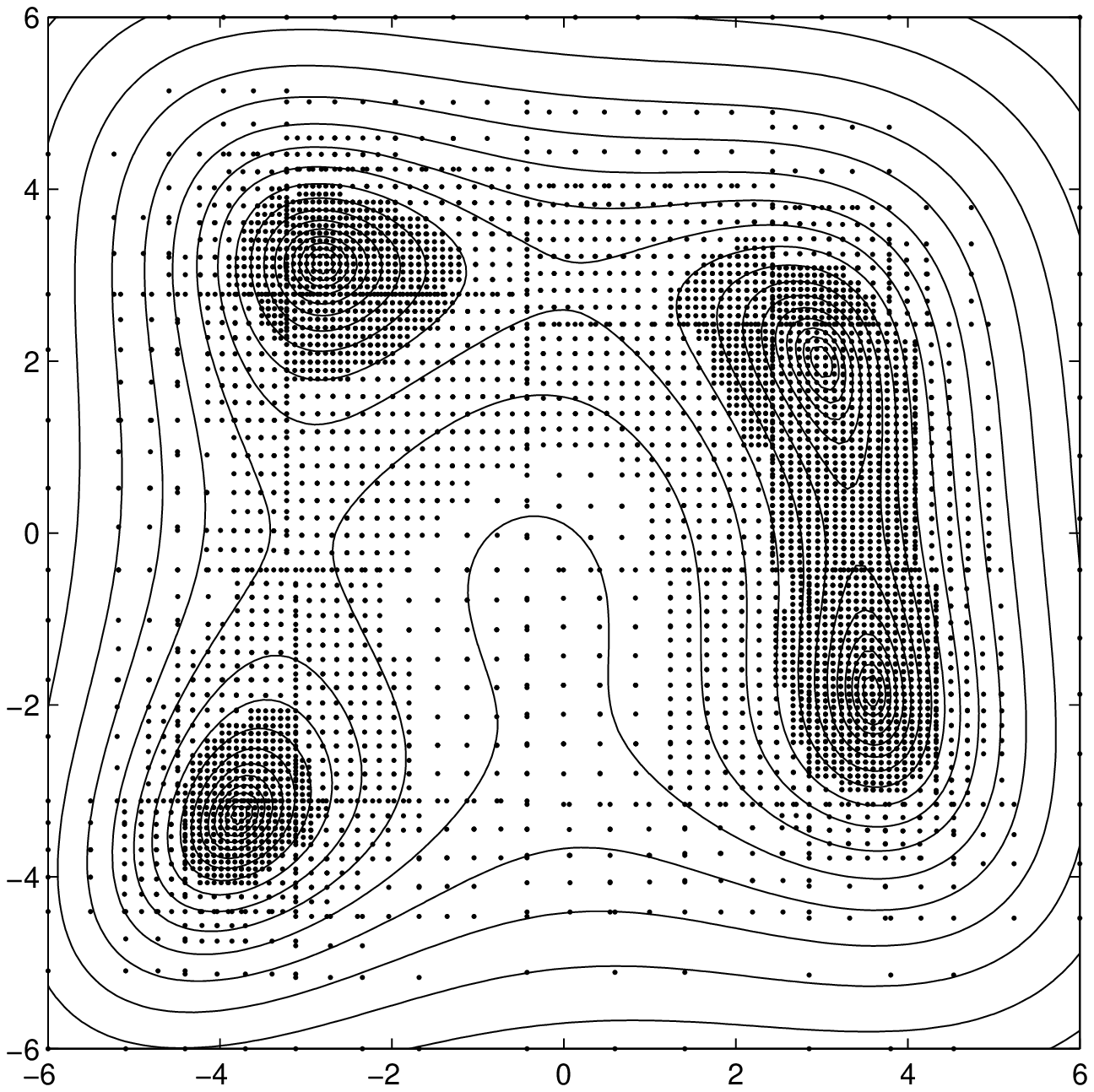,width=82mm,height=82mm,angle=0,silent=}}
\caption{Level curves of Problem~10 with the trial points
generated by strategy Partition~$2^n$ and method with global
estimate of Lipschitz constant with $r=1.1$, the number of
trials~\mbox{$=6267$}} \label{himmel11}
\end{figure}

\begin{figure}
\centerline{\psfig{file=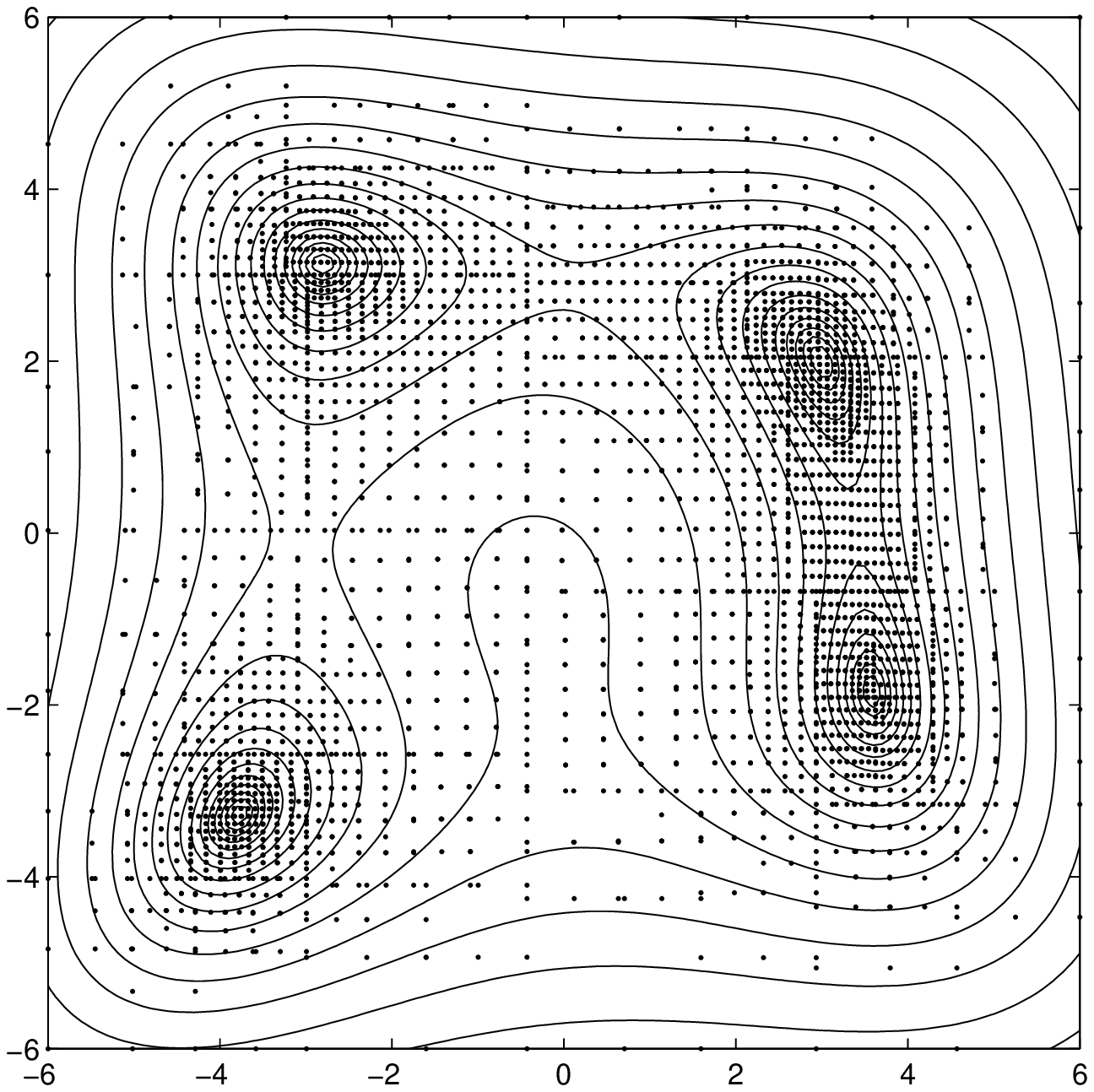,width=82mm,height=82mm,angle=0,silent=}}
\caption{Level curves of Problem~10 with the trial points
generated by strategy Bisection and method with global estimate of
Lipschitz constant with $r=1.1$, the number of
trials~\mbox{$=4384$}} \label{himmel21}
\end{figure}

\begin{figure}
\centerline{\psfig{file=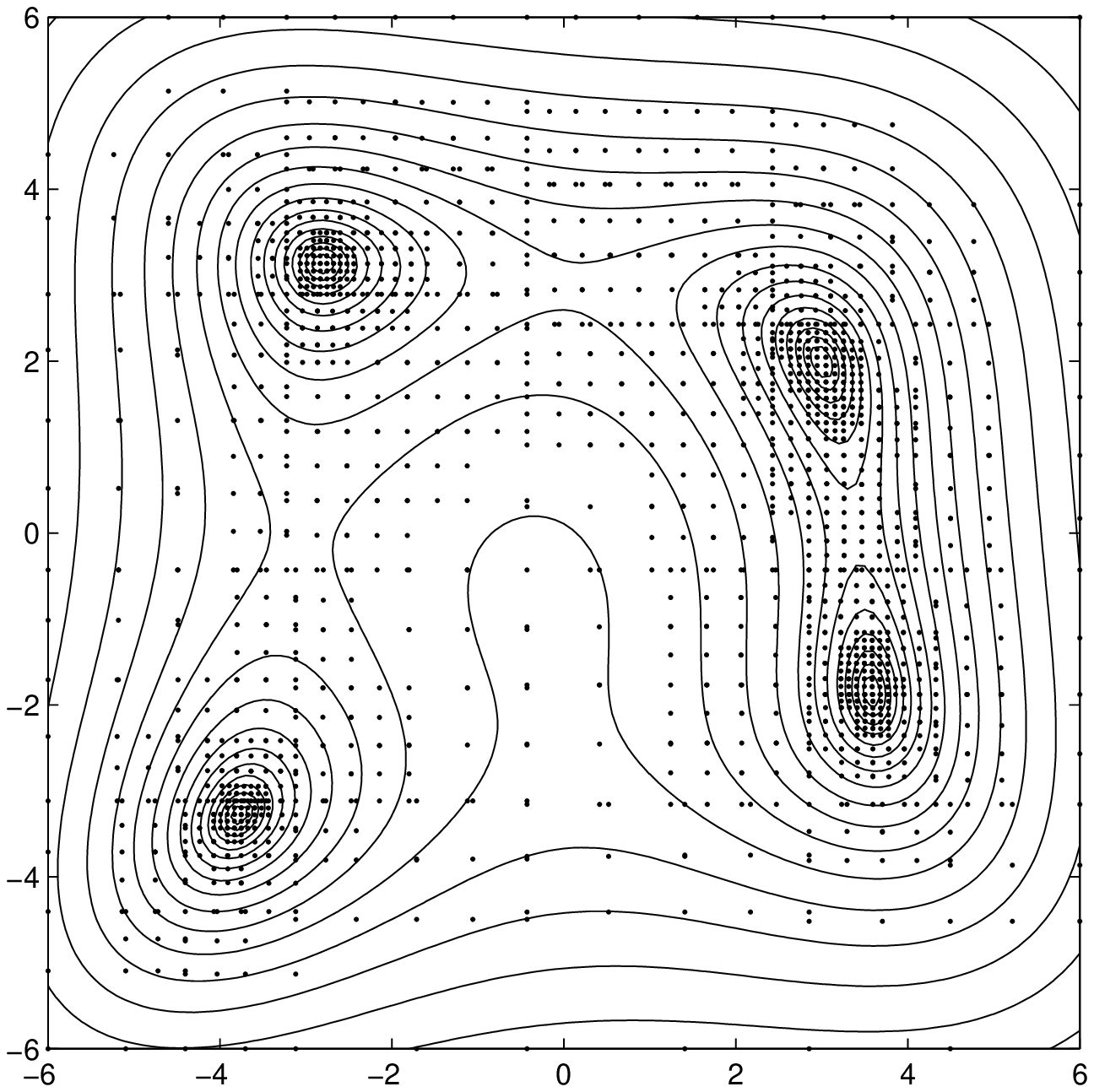,width=82mm,height=82mm,angle=0,silent=}}
\caption{Level curves of Problem~10 with the trial points
generated by strategy Partition~$2^n$ and method with local
tuning with $r=1.1$, the number of trials~\mbox{$=1797$}}
\label{himmel15}
\end{figure}

\begin{figure}
\centerline{\psfig{file=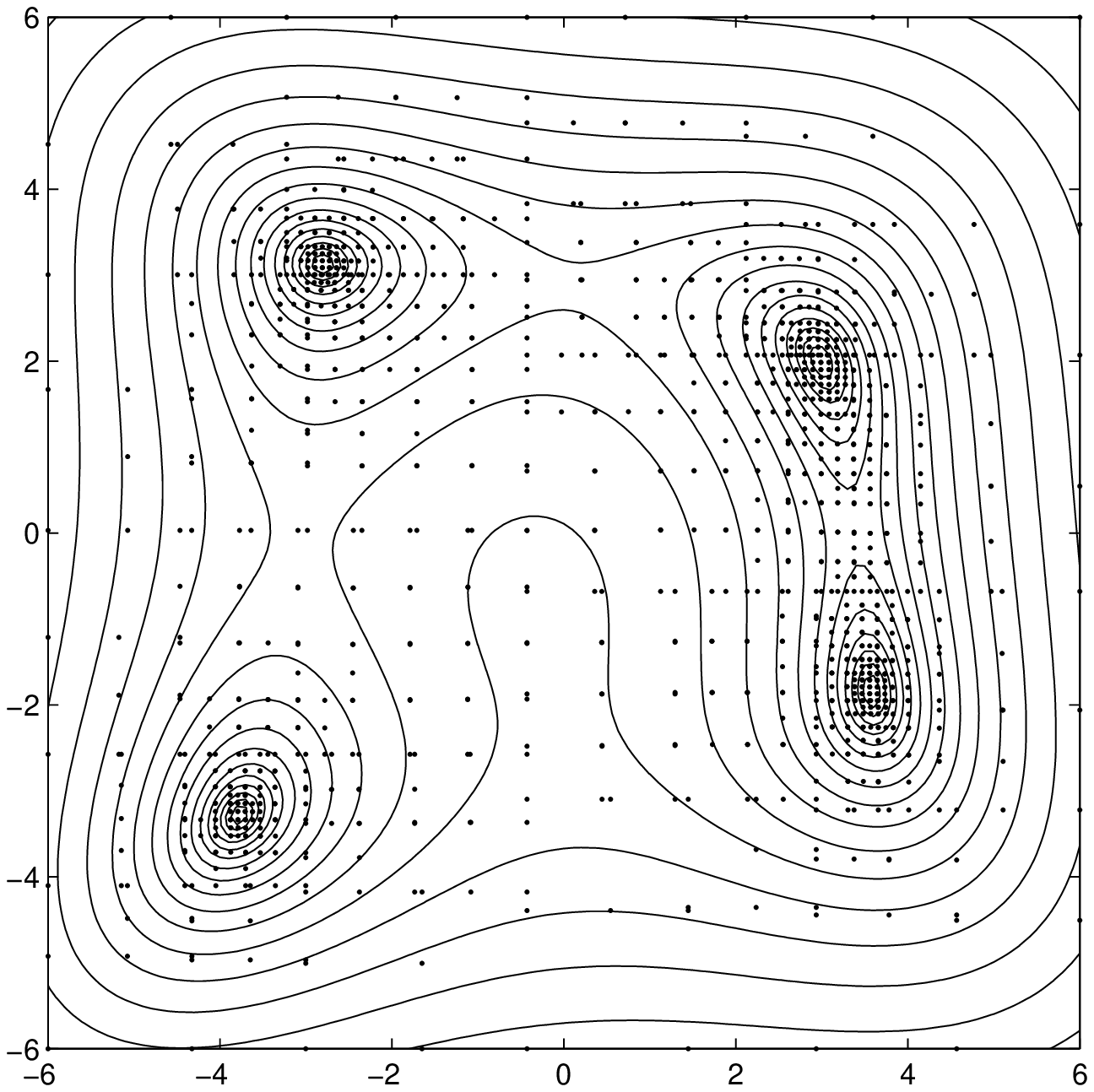,width=82mm,height=82mm,angle=0,silent=}}
\caption{Level curves of Problem~10 with the trial points
generated by strategy Bisection and method with local tuning with
$r=1.1$, the number of trials~\mbox{$=1286$}} \label{himmel25}
\end{figure}

Performance of all the methods during solving Problem~10 is
illustrated in Figs.~\ref{himmel11} -- \ref{himmel25}. Trials
points are shown by the black dots.

The new algorithm was faster than the method using the global
estimate for both strategies in all the cases. The smaller values
of the accuracy~$\varepsilon$ ensure higher values of
acceleration. For example, Table~\ref{tabProblem7} shows that the
NDAL works better when accuracy increases and the improvement is
stronger for higher values of the parameter $r$.

\begin{table}
\caption{Number of trials for Problem 7 in dependence on the
parameter $r$ and accuracy $\varepsilon$} \label{tabProblem7}
 \begin{center}
\begin{tabular}{@{\extracolsep{\fill}}|c|c|r|r|r|r|}
\hline $r$& $\varepsilon$&\multicolumn{2}{c|}{Global
Estimate}&\multicolumn{2}{c|}{Local Tuning} \\
\cline{3-6} & &Partition $2^n$&Bisection&Partition $2^n$&Bisection\\
\hline  &0.0100&6837 &4766 &5532 &3936 \\
\cline{2-6} $1.1$&0.0010 &10742 &11664 &7012 &4662 \\
\cline{2-6} &0.0001 &35697 &32218 &7367 &4694 \\
\hline
\hline &0.0100 &7267 &7374 &6422 &5128 \\
\cline{2-6} $1.3$ &0.0010 &23712 &17322 &8962 &8270 \\
\cline{2-6}  &0.0001&54397 &42584 &11862 &8582 \\
\hline
\end{tabular}
 \end{center}
\end{table}


It can be seen from the numerical experiments that the new method
with local tuning significantly outperforms the traditional
approach. In its turn, Bisection works better then Partition
$2^n$ strategy. The best combination is the new algorithm with
local tuning working with Bisection strategy.

Higher values of the parameter $r$ increase the reliability of
the methods and lead to a fast growth of the iterations number.
This happens because by increasing $r$ we uniformly augment the
estimates of the Lipschitz constants (both global and local
ones). The obtained improvement increases for higher values of
the parameter $r$.

If in the search region there exists a neighborhood of the global
solution having local Lipschitz constants smaller than the global
one (this is true, for example, for differentiable functions
having the global solution in an interior point of the search
domain), then smaller values of the accuracy $\varepsilon$ ensure
higher values of acceleration.

\end{document}